\documentclass[12pt]{amsart}
\usepackage{amsxtra,latexsym,amssymb}
\usepackage{epsfig,rotate,amsthm}
\usepackage{hyperref}

\def\qed{{\hfill $\Box$}}

\def\N{{\mathbb N}}
\def\Z{{\mathbb Z}}
\def\C{{\mathbb C}}
\def\Q{{\mathbb Q}}
\def\K{{\mathbb K}}
\def\P{{\mathcal{P}}}
\def\U{{U_{r,s}(\mathfrak sl_{\infty})}}

\def \a{{\bf a }}
\def \g{\mathfrak sl_{\infty}}
\def \B{{\mathcal{B}}}
\def \dim{{\underline{dim}}}
\theoremstyle{theorem}
\newtheorem{thm}{Theorem}[section]
\newtheorem{cor}{Corollary}[section]
\newtheorem{prop}{Proposition}[section]

\newtheorem{lem}{Lemma}[section]
\theoremstyle{definition}
\newtheorem{defn}{Definition}[section]
\theoremstyle{remark}

\begin{document}
\title[Two-parameter quantum groups] {Two-Parameter Quantum Groups and Ringel-Hall algebras of $A_{\infty}-$type}
\author[X. Tang]{Xin Tang}
\address{Department of Mathematics \& Computer Science\\
Fayetteville State University\\
1200 Murchison Road, Fayetteville, NC 28301}
\email{xtang@uncfsu.edu}
\date{\today}
\thanks{This research project is partially supported by the ISAS mini-grant at Fayetteville State University}
 \begin{abstract}
In this paper, we study the two-parameter quantum group $U_{r,s}(\mathfrak sl_{\infty})$ associated to the Lie algebra $\mathfrak sl_{\infty}$. We shall prove that 
the two-parameter quantum group $U_{r,s}(\mathfrak sl_{\infty})$ admits both a Hopf algebra 
structure and a triangular decomposition. In particular, it can be realized as the Drinfeld double 
of its certain Hopf subalgebras. We will also study a two-parameter twisted Ringel-Hall algebra $H_{r,s}(A_{\infty})$ associated to the category of all finite dimensional representations of the infinite linear quiver $A_{\infty}$. In particular, we will establish an iterated skew polynomial presentation of $H_{r,s}(A_{\infty})$ and construct a PBW basis for $H_{r,s}(A_{\infty})$.
We will establish an algebra isomorphism from $U_{r,s}^{+}(\mathfrak sl_{\infty})$ onto $H_{r,s}(A_{\infty})$. Via the theory of generic extensions in the category of finite dimensional representations of $A_{\infty}$, we shall construct several monomial bases and a bar-invariant basis for $U^{+}_{r,s}(\mathfrak sl_{\infty})$.
\end{abstract}
\maketitle                     
\section*{Introduction}
As generalizations or variations of the notation of quantum groups \cite{Dri}, several multi-parameter quantum groups have appeared in the literatures \cite{AST, CM, DP, DT, Jin, Kul, Res, Sud, Tak}. Let $\mathfrak g$ be a finite dimensional complex simple Lie algebra. Let us choose $r,s\in \C^{\ast}$ in such a way that $r,s$ are transcendental over $\Q$. The study of the two-parameter quantum group $U_{r,s}(\mathfrak g)$ has been revitalized in \cite{BGH, BKL, BW1, BW2, BW3} and the references therein. Note that the one-parameter quantum groups associated to Lie algebras ${\mathfrak gl_{\infty}}, {\mathfrak sl_{\infty}}$ of infinite ranks \cite{Kac} have been studied in the literatures \cite{DF, LS, PS1, PS2}. Similar to the case of one-parameter quantum groups, one might be interested in the constructions of the corresponding two-parameter quantum groups. 

It is the purpose of this paper to study the two-parameter quantum group $U_{r,s}(\g)$, where the Lie algebra $\g$ consists of all infinite trace-zero square matrices with only finitely many non-zero entries. We shall first define such a two-parameter quantum group and then study some of its basic properties. As usual, we will prove that the algebra $U_{r,s}(\g)$ admits a Hopf algebra structure and it is the Drinfeld double of its certain Hopf subalgebras.

To further investigate the structure of $U_{r,s}(\mathfrak sl_{\infty})$, we shall study its subalgebra $U_{r,s}^{+}(\mathfrak sl_{\infty})$ employing the approach of Ringel-Hall algebras. Note that the Ringel-Hall algebra approach has found many important applications in the study of one-parameter quantum groups \cite{Gre, Lus1, Lus2, Rei, Rin1, Rin2, Rin3, Rin4, Rin5, Xiao} and the references therein. We shall first define a two-parameter twisted Ringel-Hall algebra $H_{r,s}(A_{\infty})$ associated to the category of all finite dimensional representations of the infinite linear quiver $A_{\infty}$. Then we shall prove that the algebra $H_{r,s}(A_{\infty})$ can be presented as an iterated skew polynomial ring, and construct a PBW basis for $H_{r,s}(A_{\infty})$. Furthermore, we shall prove that $H_{r,s}(A_{\infty})$ is a direct limit of the two-parameter twisted Ringel-Hall algebras $H_{r,s}(A_{n}), n\geq 1$ associated to the finite linear quivers $A_{n}, n\geq 1$ (See \cite{Rei, T1}). 

We will establish an algebra isomorphism from $U_{r,s}^{+}(\mathfrak sl_{\infty})$ onto $H_{r,s}(A_{\infty})$. On the one hand, such an algebra isomorphism provides a generator-relation presentation of the two-parameter Ringel-Hall algebra $H_{r,s}(A_{\infty})$, which has been defined over a prescribed basis. On the other hand, via this isomorphism, we can prove that the algebra $U_{r,s}^{+}(\mathfrak sl_{\infty})$ can be presented as an iterated skew polynomial ring and it is a direct limit of $U_{r,s}^{+}(\mathfrak sl_{n+1}), n\geq 1$. As a result, we are able to construct a PBW basis for $U_{r,s}^{+}(\mathfrak sl_{\infty})$. 

To study the Borel subalgebras $U_{r,s}^{\geq 0}(\mathfrak sl_{\infty})$ of $U_{r,s}(\mathfrak sl_{\infty})$, we will study the extended two-parameter twisted Ringel-Hall algebra $\overline{H_{r,s}(A_{\infty})}$ and establish an Hopf algebra isomorphism from $U_{r,s}^{\geq 0}(\mathfrak sl_{\infty})$ onto $\overline{H_{r,s}(A_{\infty})}$. We will follow the lines in \cite{Gre, Xiao}. This result shall provide a realization of the whole two-parameter quantum group $U_{r,s}(\mathfrak sl_{\infty})$ via the double of two-parameter extended Ringel-Hall algebras. 

Note that there exists a $\Q-$algebra automorphism (which will be called the bar-automorphism) on the algebra $U_{r,s}^{+}(\mathfrak sl_{\infty})$, which exchanges $r^{\pm 1}$  and $s^{\pm 1}$ and fixes the generators $e_{i}$. Using the theory of generic extensions in the category of finite dimensional representations of $A_{\infty}$, we will construct several monomial bases for the two-parameter quantum groups following the idea used in \cite{DD1, Rei}. As an application, we will also construct a bar-invariant basis for the algebra $U^{+}_{r,s}(\mathfrak sl_{\infty})$ following \cite{Rei}.

The paper is organized as follows. In Section 1, we give the definition of $U_{r,s}(\mathfrak sl_{\infty})$ and study some of its basic properties. In Section 2, we define and study two-parameter Ringel-Hall algebra $H_{r,s}(A_{\infty})$ and establish the algebra isomorphism from $U_{r,s}^{+}(\mathfrak sl_{\infty})$ onto $H_{r,s}(A_{\infty})$.  In Section 3, we define and study the extended two-parameter Ringel-Hall algebra $\overline{H_{r,s}(A_{\infty})}$ and establish the Hopf algebra isomorphism from $U_{r,s}^{\geq 0}(\mathfrak g)$ onto $\overline{H_{r,s}(A_{\infty})}$. In Section 4, we use generic extension theory to construct some monomial bases and a bar-invariant basis for $U_{r,s}^{+}(\mathfrak sl_{\infty})$. 

\section{Definition and basic properties of the two-parameter
  quantum groups $U_{r,s}(\mathfrak sl_{\infty})$}

Let $r,s$ be two-parameters chosen from $\C^{\ast}$, such that $r,s$ are transcendental over the field $\Q$  and $r^{m}s^{n}=1$ implies $m=n=0$. Let us set $\mathcal{Z}=\Z[r^{\pm 1}, s^{\pm 1}]$ and $\mathcal{A}=\Q[r,s]_{(r-1,s-1)}$, which is the localization of $\Q[r,s]$ at the maximal ideal $(r-1, s-1)$.

Let $\mathfrak sl_{\infty}$ denote the infinite dimensional complex 
Lie algebra which consists of all trace-zero square matrices 
$(a_{ij})_{i,j\in \N}$ with only finitely many non-zero entries. 
The one-parameter quantum groups $U_{q}(\mathfrak sl_{\infty})$ 
associated to $\mathfrak sl_{\infty}$ were studied by various people 
in the references \cite{DF, LS, PS1, PS2}. Following a similar idea in \cite{BW1, Tak}, 
we will introduce a class of two-parameter quantum groups $U_{r,s}(\mathfrak \g)$ 
associated to the Lie algebra $\mathfrak sl_{\infty}$.

It is well known that one can also define roots for the Lie algebra 
$\mathfrak sl_{\infty}$ as in the finite dimensional case of ${\mathfrak
g}={\mathfrak sl_{n+1}}$. In particular, all the simple roots of $\mathfrak \g$ can be 
denoted as $\alpha_{i}, i\in I=\N$. Accordingly, all the positive roots of $\mathfrak \g$ 
are exactly given as $\alpha_{ij}\colon =\sum _{k=i}^{j}\alpha_{k}$ for $i\leq j \in \N$. 

Let $C=(c_{ij})_{i,j\in \N}$ denote the infinite Cartan matrix corresponding to the Lie algebra $\mathfrak sl_{\infty}$. Then, we have the following
\[
c_{ii}=2, \, c_{ij}=-1\,\text{for}\, |i-j|=1, \, c_{ij}=0 \, \text{for}\,
|i-j|>1.
\]
Let $\Q(r,s)$ denote the function field in two variables $r,s$ over
the field $\Q$ of all rational numbers. Let $\mathcal{Q}$ denote the root 
lattice generated by $\alpha_{i}, i\in \N$. Then we can define a bilinear 
form $\langle -,-\rangle$ on the root lattice $\mathcal{Q}\cong \Z^{\oplus \N}$ 
as follows
\[
\langle i,j  \rangle \colon =\langle \alpha_{i},\alpha_{j} \rangle
=\left
\{\begin{array}{cc}a_{ij}, & {\rm if}\, i<j,\\
1,& {\rm if} \,i=j,\\
0, & {\rm if } \,i>j.
\end{array}
\right.
\]

\begin{defn} 
The two-parameter quantum group $U_{r,s}(\mathfrak sl_{\infty})$ is 
defined to be the $\Q(r,s)-$algebra generated by $e_{i},f_{i}, w_{i}^{\pm 1}, w_{i}^{\prime\pm 1}, i\in \N$ subject to the 
following relations
\begin{eqnarray*}
w_{i}^{\pm 1}w_{j}^{\pm 1}=w_{j}^{\pm 1}w_{i}^{\pm 1},
\quad 
w_{i}^{\prime\pm 1}w_{j}^{\prime \pm 1}=w_{j}^{\prime \pm 1}
w_{i}^{\prime \pm 1},\\
w_{i}^{\pm 1}w_{j}^{\prime \pm 1}= w_{j}^{\prime \pm 1}w_{i}^{\pm 1},
\quad 
w_{i}^{\pm 1}w_{i}^{\mp 1}=1= w_{i}^{\prime\pm 1}w_{i}^{\prime
  \mp1},\\
w_{i}e_{j}=r^{\langle j,i\rangle }s^{-<i,j>}e_{j}w_{i},
\quad 
w^{\prime}_{i}e_{j}=r^{-\langle i,j\rangle }s^{\langle j,i\rangle }e_{j}w^{\prime}_{i},\\
w_{i}f_{j}=r^{-\langle j,i\rangle }s^{\langle i,j\rangle}f_{j}e_{i},\quad w^{\prime}_{i}f_{j}=r^{\langle i,j\rangle}s^{-\langle j,i\rangle }f_{j}w^{\prime}_{i},\\
e_{i}f_{j}-f_{j}e_{i}=\delta_{i,j}\frac{w_{i}-w^{\prime}_{i}}{r_{i}-s_{i}},\\
e_{i}e_{j}-e_{j}e_{i}=f_{i}f_{j}-f_{j}f_{i}=0 \, \text{for}\, |i- j|>1,\\
e_{i}^{2}e_{i+1}-(r+s)e_{i}e_{i+1}e_{i}+rse_{i+1}e_{i}^{2}=0,\\
e_{i}e_{i+1}^{2}-(r+s)e_{i+1}e_{i}e_{i+1}+rse_{i+1}^{2}e_{i}=0,\\
f_{i}^{2}f_{i+1}-(r^{-1}+s^{-1})f_{i}f_{i+1}f_{i}+(rs)^{-1}f_{i+1}f_{i}^{2}=0,\\
f_{i}f_{i+1}^{2}-(r^{-1}+s^{-1})f_{i+1}f_{i}f_{i+1}+(rs)^{-1}f_{i+1}^{2}f_{i}=0.\\
\end{eqnarray*}
\end{defn}

First of all, we have the following obvious proposition concerning a Hopf algebra structure of the algebra $U_{r,s}(\mathfrak sl_{\infty})$.
\begin{prop}
The algebra $U_{r,s}(\mathfrak g)$ is a Hopf algebra with the comultiplication, counit and 
antipode defined as follows
\begin{eqnarray*}
\Delta(w_{i}^{\pm 1})=w_{i}^{\pm 1}\otimes w_{i}^{\pm 1},\quad 
\Delta(w_{i}^{\prime\pm 1})=w_{i}^{\prime\pm 1}\otimes w_{i}^{\prime\pm 1},\\
\Delta(e_{i})=e_{i}\otimes 1+ w_{i}\otimes e_{i}, \quad 
\Delta(f_{i})=1\otimes f_{i} + f_{i}\otimes w^{\prime}_{i},\\
\epsilon ( w_{i}^{\pm 1})=\epsilon (w_{i}^{\prime\pm 1})=1, \quad 
\epsilon (e_{i})=\epsilon(f_{i})=0,\\
S(w_{i}^{\pm 1})=w_{i}^{\mp 1},\quad 
S(w_{i}^{\prime\pm 1})=w_{i}^{\prime \mp 1},\\
S(e_{i})=-w_{i}^{-1}e_{i},\quad S(f_{i})=-f_{i}w_{i}^{\prime-1}.
\end{eqnarray*}
\end{prop}
{\bf Proof:} The proof is reduced to the finite case where $\mathfrak{g}={\mathfrak sl_{n+1}}$, whose proof can be found in \cite{BW1}. And we will not repeat the details here.
\qed

Let $U_{r,s}^{+}(\mathfrak sl_{\infty})$ (resp. $U_{r,s}^{-}(\mathfrak sl_{\infty})$) denote the subalgebra of $U_{r,s}(\mathfrak sl_{\infty})$ generated by $e_{i}, i\in \N$ (resp. by $f_{i}, i\in \N$). Let $U_{r,s}^{0}(\mathfrak sl_{\infty})$ denote the subalgebra of
$U_{r,s}(\mathfrak sl_{\infty})$ generated by $w_{i}^{\pm 1},
w_{i}^{\prime \pm 1}, i\in \N$. Then we shall have the following triangular decomposition of $U_{r,s}(\mathfrak sl_{\infty})$.
\begin{prop}
The algebra $U_{r,s}(\mathfrak sl_{\infty})$ has a triangular decomposition
\[
U_{r,s}(\mathfrak sl_{\infty})\cong U^{-}_{r,s}(\mathfrak sl_{\infty})\otimes U^{0}_{r,s}(\mathfrak
sl_{\infty})\otimes 
U^{+}_{r,s}(\mathfrak sl_{\infty}).
\]
\end{prop}
{\bf Proof:} Once again, we can repeat the proof used in the case of $U_{r,s}(\mathfrak sl_{n+1})$. We refer the reader to \cite{BW1} for more details.
\qed 

Let us denote by $\Z^{\oplus \N}$ the free abelian group of rank $|\N|$ with a
basis denoted by $z_{1}, z_{2}, \cdots, z_{n},\cdots$. Given any element 
${\bf a} \in \Z^{\oplus \N}$, say ${\bf a}=\sum a_{i}z_{i}$, we set $|{\bf a}|
=\sum a_{i}$. Note that algebra $U^{+}_{r,s}(\mathfrak \g)$ (resp. 
$U^{-}_{r,s}(\mathfrak sl_{\infty})$) is a $\Z^{\oplus \N}-$graded algebra 
by assigning to the generator $e_{i}$ (resp. $f_{i}$) the degree $z_{i}$. Given
${\bf a}\in \Z^{\oplus \N}$, we denote by $U_{r,s}^{\pm}(\mathfrak
sl_{\infty})_{\bf a}$ the set of homogeneous elements of degree ${\bf a}$ in
$U^{\pm}_{r,s}(\mathfrak sl_{\infty})$.
\begin{prop} We have the following decomposition
\[
U_{r,s}^{+}(\mathfrak \g)=\bigoplus_{\bf a}U_{r,s}^{+}(\mathfrak \g)_{\bf
a},\quad U_{r,s}^{-}(\mathfrak \g)=\bigoplus_{\bf a}
U_{r,s}^{-}(\mathfrak \g)_{\bf a}.
\] 
\end{prop}
\qed

Let us define $U_{v,v^{-1}}(\mathfrak sl_{\infty})$ to be the specialization of $U_{r,s}(\mathfrak sl_{\infty})$ for $r=v=s^{-1}$. Then we shall have the following similar result as \cite{BW1}.
\begin{prop}
Assume there exists an isomorphism of Hopf algebras
\[
\phi \colon U_{r,s}(\mathfrak sl_{\infty})\longrightarrow U_{v,v^{-1}}(\mathfrak sl_{\infty})
\]
for some $v$. Then $r=v$ and $s=v^{-1}$.
\end{prop}
\qed

\subsection{A Drinfeld double realization of $U_{r,s}(\mathfrak \g)$}
In this subsection, we show that the two-parameter quantum group $U_{r,s}(\frak \g)$ can be realized as the Drinfeld double of its certain Hopf subalgebras. To proceed, we need to recall a couple of standard definitions for the Hopf pairing and the Drinfeld double of Hopf algebras. For more details about these concepts, we refer the reader to the references \cite{BW1,DT}.

\begin{defn} A Hopf pairing of two Hopf algebras $H^{\prime}$
  and $H$ is a bilinear form $(,)\colon H^{\prime}  \times
  H\longrightarrow \K$ such that 

\begin{enumerate}
\item  $(1,h)=\epsilon_{H}(h),$\\

\item $(h^{\prime}, 1)=\epsilon_{H^{\prime}}(h^{\prime}),$\\

\item $(h^{\prime},hk)=(\Delta_{H^{\prime}}(h^{\prime}),h\otimes k)=\sum (h^{\prime}_{(1)},h)(h^{\prime}_{(2)},k),$\\

\item $(h^{\prime}k^{\prime}, h)=(h^{\prime}\otimes k^{\prime}, \Delta_{\prime}(h))=\sum (h^{\prime}, h_{(1)})(k^{\prime},h_{(2)}),$
\end{enumerate}
for all $h,k\in H^{\prime}$,$ h^{\prime}, k^{\prime}\in  H^{\prime}$, where $\epsilon_{H}, \epsilon_{H^{\prime}}$ denote the counits of $H, H^{\prime}$ respectively, and $\Delta_{H},\Delta_{H^{\prime}}$ denote their comultiplications.
\end{defn}

It is obvious that 
\[
(S_{H^{\prime}}(h^{\prime}), h)=(h^{\prime},S_{H}(h))
\]
for all $h\in H$ and $h^{\prime} \in H^{\prime}$, where $S_{H^{\prime}}$ and $S_{H}$ denote the respective antipodes of $ H$ and $H^{\prime}$.

Let $U_{r,s}^{\geq 0}(\mathfrak \g)$ (resp. $U_{r,s}^{\leq 0}(\mathfrak \g)$) be the Hopf subalgebra of $U_{r,s}(\mathfrak \g)$ generated by $e_{i}, w_{i}^{\pm
  1}$ (resp. $f_{i},w_{i}^{\prime \pm 1}$). Assume that $B=U_{r,s}^{\geq 0}(\mathfrak \g)$ 
and $(B^{\prime})^{coop}$ is the Hopf algebra generated by $f_{j}, (w_{j}^{\prime})^{\pm 1}$ with the opposite coproduct to $U^{\leq 0}(\mathfrak sl_{\infty})$. Using the same proof  in the case of $\mathfrak sl_{n+1}$ \cite{BW1}, we shall have the following result
\begin{lem}
There exists a unique Hopf pairing $B$ and $B^{\prime}$ such that 
\begin{eqnarray*}
(f_{i}, e_{j})=\frac{\delta_{i,j}}{s-r}\\
(w_{i}^{\prime},w_{j})=r^{<e_{i},e_{j}>}s^{-<e_{j}, e_{i}>},
\end{eqnarray*}
and the pairing takes the zero value on all other pairs of generators. Moreover, we have
$(S(a),S(b))=(a,b)$ for $a\in B^{\prime}, b\in B$.
\end{lem}
\qed 

Therefore, we have the following similar result as in \cite{BW1}.
\begin{thm}
$U_{r,s}(\mathfrak \g)$ can be realized as a Drinfeld double of Hopf
subalgebras $B=U^{\geq 0}_{r,s}(\mathfrak \g)$ and $(B^{\prime})^{coop}=U^{\leq 0}_{r,s}(\mathfrak \g)$, that is,
\[
U_{r,s}(\mathfrak \g)\cong D(B, (B^{\prime})^{coop}).
\]
\end{thm}
{\bf Proof:} First of all, let us define a linear map: $\phi \colon D(B,
(B^{\prime})^{coop})\longrightarrow U_{r,s}(\mathfrak sl_{\infty})$ as follows
\[
\phi(\hat{\omega}_{i}^{\pm 1})=\omega_{i}^{\pm 1}, \quad
\phi((\hat{\omega}_{i}^{\prime})^{\pm 1})=(\omega_{i}^{\prime})^{\pm 1}
\]
\[
\phi(\hat{e}_{i})=e_{i},\quad \phi_{i}(\hat{f}_{i})=f_{i}.
\]
We need to show that this mapping is a Hopf algebra automorphism. Obviously, we can still 
employ the proof used in \cite{BW1} for the finite case $g=\mathfrak sl_{n+1}$ and we will not repeat the detail here.
\qed

\subsection{An integral form of the two-parameter quantum group $\U$} In addition, we can consider an integral form of the two-parameter quantum group $\U$ and its subalgebras following \cite{Lus1}. For any $l\geq 1$, let us set the following
\[
[l]=\frac{r^{l}-s^{l}}{r-s}, \, [l]^{!}=[1][2]\cdots [l].
\]
Let us define $e_{i}^{(l)}=\frac{e_{i}^{l}}{[l]^{!}}, f_{i}^{(l)}=\frac{f_{i}^{l}}{[l]^{!}}$. 
We define a $\mathcal{Z}-$subalgebra $U_{r,s}(\mathfrak sl_{\infty})_{\mathcal{Z}}$ of $\U$ which is generated by the elements $e_{i}^{(l)}, f_{i}^{(l)}, w_{i}^{\pm 1}, w_{i}^{\prime \pm 1}$ for $i\in I$. Similarly, we can define the integral form of $U_{r,s}^{+}(\mathfrak sl_{\infty})$ and $U_{r,s}^{-}(\mathfrak sl_{\infty})$. It is easy to see that we have the  following
\[
\U\cong U_{r,s}(\mathfrak sl_{\infty})_{\mathcal{Z}}\otimes_{\mathcal{Z}} \Q(r,s)
\]
and 
\[
U^{\pm}_{r,s}(\mathfrak sl_{\infty})\cong U_{r,s}^{\pm}(\mathfrak sl_{\infty})_{\mathcal{Z}} \otimes_{\mathcal{Z}}\Q(r,s).
\]
In particular, $\U$ (resp. $U_{r,s}^{\pm }(\mathfrak sl_{\infty})$) is a free $\mathcal{Z}-$algebra.

\section{Two-parameter Ringel-Hall algebras $H_{r,s}(A_{\infty})$}
To study the two-parameter quantum group $U_{r,s}(\mathfrak \g)$, it is helpful to study its subalgebra $U_{r,s}^{+}(\mathfrak \g)$. We shall study this algebra in terms of two-parameter Ringel-Hall algebra associated to the infinite linear quiver. We will define and study a two-parameter Ringel-Hall algebra $H_{r,s}(A_{\infty})$ associated to the category of finite 
dimensional representations of the infinite quiver
\[
A_{\infty}\colon \stackrel{1}{\bullet} \longrightarrow \stackrel{2}{\bullet} \longrightarrow\stackrel{3}{\bullet} \cdots \stackrel{n-1}{\bullet} \longrightarrow\stackrel{n}{\bullet} \longrightarrow \cdots.
\]

For $n\geq 1$, let $A_{n}$ denote the finite quiver 
\[
A_{n}\colon \stackrel{1}{\bullet} \longrightarrow \stackrel{2}{\bullet} \longrightarrow\stackrel{3}{\bullet} \cdots\stackrel{n-1}{\bullet} \longrightarrow \stackrel{n}{\bullet}
\]
with $n$ vertices. Let us fix $k$ to be a finite field and let $\Lambda_{n}$ denote the path 
algebra of the finite linear quiver $A_{n}$ over $k$. Then $\Lambda_{n}$ is a finite dimensional 
hereditary algebra of finite-representation type. Note that the category of
finite dimensional representations of the quiver $A_{n}$ is equivalent to the
category of finite dimensional $\Lambda_{n}-$modules. We will denote
this category by $A_{n}-$mod. Let us set $q=|k|$ the cardinality of $k$, 
and choose $v$ to be a number such that $v^{2}=q$. We know that 
$\Lambda_{n}$ is finitary in the sense that the cardinality of 
the extension group $Ext^{1}(S,S^{\prime})$ is finite for any 
two simple $\Lambda_{n}-$modules $S, S^{\prime}$.

Let us denote by $A_{\infty}-$mod, the category of all finite dimensional representations of the quiver $A_{\infty}$. Note that the category $A_{\infty}$-mod has been investigated by Hou 
and Ye in \cite{HY}, where they have explicitly described all finite dimensional indecomposable representations of $A_{\infty}$ and studied the one-parameter non-twisted generic Ringel-Hall algebra $H_{q}(A_{\infty})$. Let $S_{i}$ be the simple representation associated to the vertex 
$i$ of the quiver $A_{\infty}$ and let $M_{ij}$ denote the indecomposable representation of $A_{\infty}$ with a top $S_{i}$ and length $j-i+1$. It is easy to see that there is a one-one to correspondence between the set of isoclasses of finite dimensional indecomposable representations $M_{ij}$ of the quiver $A_{\infty}$ and the set of positive roots $\alpha_{ij}$ 
for the Lie algebra $\mathfrak \g$. 

Concerning the relationship between the categories $A_{n}-$mod and $A_{\infty}-$mod, we now recall the following result from \cite{HY}.
\begin{thm} 
({\bf Theorem 1.1} in \cite{HY}) The category $A_{n}-$mod can be regarded as a fully 
faithful and extension closed subcategory of $A_{\infty}-$mod 
and $A_{m}-$mod for $m\geq n$.
\end{thm}
\qed

Based on the above theorem, we know that the extension group between any two finite dimensional representations $M, N$ of $A_{\infty}$ can be calculated via regarding $M,N$ as the 
representations of a certain finite quiver $A_{m}$. Therefore, the number of extensions between $M, N$ is still depicted by the evaluation of  the Hall polynomial at $q$, the cardinality of the base field. Recall that the two-parameter Ringel-Hall algebra $H_{r,s}(A_{n}), n\geq 1$ associated to the category $A_{n}-$mod has been studied in \cite{Rei, T1}. In particular, one knows that $H_{r,s}(A_{n})$ can be presented as an iterated skew polynomial ring and its prime ideals are completely prime. A PBW basis has also been constructed for $H_{r,s}(A_{n})$ in \cite{T1} as well. Note that this approach is plausible because of the existence of Hall polynomials in the category $A_{\infty}-$mod. Indeed, we will be looking at a limit version $H_{r,s}(A_{\infty})$ of the two-parameter Ringel-Hall algebras $H_{r,s}(A_{n}), n \geq 1$.

\subsection{Two-parameter Ringel-Hall algebra $H_{r,s}(A_{\infty})$} We will denote by $\mathcal{P}$ the set of isomorphism classes of finite dimensional  representations of the infinite quiver $A_{\infty}$. Let us define the subset
\[
\mathcal{P}_{1}=\mathcal{P}-{0}
\]
where $0$ denotes the subset of $\mathcal{P}$ consisting of the 
only isomorphism class of the zero representation. For any 
$\alpha \in \mathcal{P}$, we choose a representation $u_{\alpha}$ 
corresponding to $\alpha$. We denote by $a_{\alpha}$ the order of 
the automorphism group $Aut(u_{\alpha})$. It is easy to see that the 
number $a_{\alpha}$ is independent of the choices of the
representatives $u_{\alpha}$ for any $\alpha \in \mathcal{P}$. 

For any given three representatives $u_{\alpha}, u_{\beta}, u_{\gamma}$ 
of the elements $\alpha, \beta, \gamma \in \mathcal{P}$ respectively, we 
denote by $g_{\alpha\beta}^{\gamma}$ the number of submodules $N$ of
$u_{\gamma}$ satisfying the conditions: $N\cong u_{\beta}$ and $u_{\gamma}/N \cong u_{\alpha}$. 

Note that it does not make sense to define $Ext^{1}(M, N)$ for any two given representations $M, N$ of the infinite quiver $A_{\infty}$. Let us denote by $\hat{E}_{A_{\infty}}(M, N)$ the set of all short exact sequences $0\longrightarrow N \longrightarrow E\longrightarrow M\longrightarrow 0$. We say two such short exact sequences $0\longrightarrow N \longrightarrow E_{1}\longrightarrow M\longrightarrow 0$ and $0\longrightarrow N \longrightarrow E_{2}\longrightarrow M\longrightarrow 0$ are equivalent if there exists a homomorphism $\phi\colon E_{1}\longrightarrow E_{2}$ making the diagram commute. We denote by $E_{A_{\infty}}(M, N)$ the set of all equivalence classes of $\hat{E}_{A_{\infty}}$ with respect to this equivalence relation.  For any given $M, N\in A_{\infty}-mod$, according to {\bf Theorem 1.2} in \cite{HY}, we can choose some $m\geq 1$ such that there exists a bijection between $E_{A_{\infty}}(M, N)$ and $Ext^{1}_{A_{m}-mod}(M, N)$. If no confusion arises, we will still write $E_{A_{\infty}}(M, N)$ as $Ext^{1}(M, N)$. 

For any given $M, N\in A_{\infty}-mod$, we define the following notation
\[
\langle M,N\rangle =dim_{k}Hom(M,N)-dim_{k}Ext^{1}(M,N).
\]

Once the representations $M, N$ are chosen, we can always restrict to a subcategory $A_{n}-mod$. Since the algebra $\Lambda_{n}$ is hereditary for any $n\in \N$, it is easy to see 
that for any representations $M, N \in A_{n}-mod$, the value of $\langle M, N\rangle$ 
solely depends on the dimension vectors $\dim M,\, \dim N$ of the $A_{n}-$modules 
$M$ and $N$. 

Now for any given elements $\alpha, \beta \in \mathcal{P}$, we can define the following notation
\[
\langle \alpha,\beta \rangle=\langle u_{\alpha},u_{\beta}\rangle
\]
where $u_{\alpha}, u_{\beta}$ are any chosen representatives 
of $\alpha, \beta$ respectively. It is easy to see that $\langle-,-\rangle$ is a bilinear form. 

It is well known that in the category $A_{n}-$mod, there exists a
symmetry between the objects of $A_{n}-$mod. This symmetry is
described by Green's formula \cite{Gre}. In fact, one can also prove that 
Green's formula holds for the objects in the category $A_{\infty}-$mod. Namely, 
we have the following result.
\begin{thm}
Let $\alpha, \beta, \alpha^{\prime}, \beta^{\prime}\in 
\mathcal{P}$, then we have
\begin{eqnarray*}
a_{\alpha}a_{\beta}a_{\alpha^{\prime}}a_{\beta^{\prime}}\sum_{\lambda
  \in \mathcal{P}}
  g^{\lambda}_{\alpha,\beta}g^{\lambda}_{\alpha^{\prime}
  \beta^{\prime}}a_{\lambda}^{-1}=\sum_{\rho, \sigma,\sigma^{\prime},
  \tau \in \mathcal{P}}
  &\frac{|Ext^{1}(u_{\rho},u_{\tau})|}{|Hom(u_{\rho},u_{\tau})|}&
  g^{\alpha}_{\rho \sigma}g^{\alpha^{\prime}}_{\rho
  \sigma^{\prime}}g^{\beta}_{\sigma^{\prime}\tau}g^{\beta^{\prime}}_{\sigma
  \tau}\\&&a_{\rho}a_{\sigma}a_{\sigma^{\prime}}a_{\tau^{\prime}}.
\end{eqnarray*}
\end{thm}
{\bf Proof:} Since all representations involved in the formula are finite dimensional
representations of $A_{\infty}$, we can choose some positive integer $m$ such that 
$\alpha, \beta, \alpha^{\prime}, \beta^{\prime}$ and $\lambda$ can actually 
be regarded as objects in the subcategory $A_{m}-$mod instead. Note that 
Green's formula holds within the subcategory $A_{m}-$mod. Since the category 
$A_{m}-$mod is a fully faithful and extension closed subcategory of $A_{\infty}-$mod, we know that Green's formula holds in $A_{\infty}-$mod.
\qed

Let $H_{r,s}(A_{n})$ denote the two-parameter Ringel-Hall algebra
associated to the category $A_{n}-$mod as defined in \cite{Rei}. 
In \cite{Rei}, Reineke has proved that the two-parameter Ringel-Hall 
algebra $H_{r,s}(A_{n})$ is isomorphic to the algebra $U^{+}_{r,s}(\mathfrak sl{n+1})$. In the rest of this section, we will show that a limit version of this statement is still true.

Note that there exist Hall polynomials $F_{M,N}^{L}(x)$ for $M,N,L \in A_{n}-mod$ such that $g^{L}_{M,N}=F_{M,N}^{L}(q)$, where $q$ is the cardinality of the base field $k$. For the existence and calculation of Hall polynomials in $A_{n}-mod$, we refer the reader 
to the references \cite{Rin2, Rin3}. Since each $A_{n}-$mod is a fully faithful and extension closed subcategory of $A_{\infty}-$mod, the Hall polynomials exists for objects in $A_{\infty}-$mod, which leads to the definition of two-parameter Ringel-Hall algebra $H_{r,s}(A_{\infty})$ below.

Now let us define $H_{r,s}(A_{\infty})$ to be the free $\Q(r,s)-$module 
generated by the set $\{u_{\alpha}\mid \alpha \in \mathcal{P}\}$. Moreover, we define a multiplication on the 
free $\Q(r,s)-$module $H_{r,s}(A_{\infty})$ as follows
\[
u_{\alpha} u_{\beta} =\sum_{\lambda \in \mathcal{P}}
s^{-\langle \alpha ,\beta \rangle}F^{u_{\lambda} }_{u_{\alpha} u_{\beta}}(rs^{-1})u_{\lambda},\quad \text{for
  any}\, \alpha, \beta \in \mathcal{P}.
\]

It is easy to see that we have the following result.
\begin{thm}
The free $\Q(r,s)-$module $H_{r,s}(A_{\infty})$ is an associative 
$\Q(r,s)-$algebra under the above defined multiplication. In
particular, the algebra $H_{r,s}(A_{n})$ can be regarded as a 
subalgebra of $H_{r,s}(A_{\infty})$ and $H_{r,s}(A_{m})$ 
for $m\geq n$. In particular, we have 
\[
H_{r,s}(A_{\infty})=\lim_{n\mapsto \infty} H_{r,s}(A_{n}).
\]
\end{thm}
{\bf Proof:} It is straightforward to verify that $H_{r,s}(A_{\infty})$ is an associative algebra under the above defined multiplication. Once again, we can reduce the proof to the finite 
case thanks to {\bf Theorem 1.1} in \cite{HY}.  Since each $A_{n}-$mod can be regarded as a fully faithful and extension closed subcategory of $A_{\infty}-$mod and $A_{m}-$mod when $m\geq n$, the algebra  $H_{r,s}(A_{n})$ can be regarded as a subgroup of the algebras 
$H_{r,s}(A_{\infty})$ and $H_{r,s}(A_{m})$. Furthermore, one notices 
that the multiplication of $H_{r,s}(A_{n})$ is the restriction of the multiplications 
of $H_{r,s}(A_{\infty})$ and $H_{r,s}(A_{m})$. Therefore, the algebra 
$H_{r,s}(A_{n})$ can be regarded as a subalgebra of $H_{r,s}(A_{\infty})$ 
and $H_{r,s}(A_{m})$ for $m\geq n$ as desired. Furthermore, each element of 
$H_{r,s}(A_{\infty})$ can be regarded as an element of a certain subalgebra 
$H_{r,s}(A_{m})$. Thus we shall have $H_{r,s}(A_{\infty})=\lim_{n\mapsto \infty} H_{r,s}(A_{n})$ as desired.
\qed

\subsection{Basic properties of $H_{r,s}(A_{\infty})$}
Since the category $A_{\infty}-$mod can be regarded the direct limit of its fully faithful and extension closed subcategories $A_{n}-$mod with $n\geq 1$, any two objects $M, N \in A_{\infty}-mod$ can be regarded as objects in a certain subcategory $A_{m}-$mod. 
Thus the extension between any such two objects can be handled in this subcategory $A_{n}-$mod as well. As a result, it is no surprise that the algebra $H_{r,s}(A_{\infty})$ 
shares many similar ring-theoretic properties with its subalgebras $H_{r,s}(A_{n})$. In this subsection, we will establish some similar results for $H_{r,s}(A_{\infty})$ without giving detailed proofs. The reader shall be reminded that all the proofs can be
reconstructed the same way as in the case of a certain subalgebra 
$H_{r,s}(A_{m})$. And we refer the curious reader to \cite{T1} for the details.

First of all, let us fix more notations. For any given $\alpha\in \P$, we will 
choose an element $u_{\alpha} \in H_{r,s}(\Lambda)$. We denote by 
$\epsilon(\alpha)$ the $k-$dimension of the endomorphism ring of the
representative $u_{\alpha}$ associated to $\alpha$. For any given finite 
dimensional representation $M$ of the infinite quiver $A_{\infty}$, we will denote the 
isomorphism class of $M$ by $[M]$ and the dimension vector of 
$M$ by $\dim M$, which is an element of the Grothendieck group 
$K_{0}(A_{\infty})$ of the category $A_{\infty}-$mod.

Recall that there is a one-to-one correspondence between the set of
all positive roots for the Lie algebra $\mathfrak \g$ and the set of 
isoclasses of finite dimensional indecomposable representations of 
$A_{\infty}$. Let $\a \in \Phi^{+}$ be any positive root, we shall 
denote by $M(\a)$ the indecomposable representation corresponding to
$\a$. For any given map $\alpha \colon \Phi^{+}\longrightarrow \N_{0}$ with finite support, let us 
set the following 
\[
M(\alpha)=M_{\Lambda}(\alpha)=\bigoplus_{\bf a\in \Phi^{+}} \alpha(\bf a)M(\bf a).
\]
Then it is easy to see there is a one-to-one correspondence between the set $\mathcal{P}$ of isomorphism classes of all 
finite dimensional representations of the infinite quiver $A_{\infty}$ and the set of all maps $\alpha \colon 
\Phi^{+}\longrightarrow \N_{0}$ with finite supports. From now on, we 
will not distinguish an element $\alpha \in \P$ from the corresponding 
map associated to $\alpha$, and we may denote both of them 
by $\alpha$ if no confusion arises. 

For any given $\alpha \in\mathcal{P}$, let us set ${\bf dim}
\alpha=\sum_{\bf a \in \Phi^{+}} \alpha(\bf a)\a$. Then we shall 
have following
\[
\dim M(\alpha)={\bf dim } \alpha.
\]

For any given $\alpha \in \P$, we will denote by $dim (\alpha)=dim(u_{\alpha})$ 
the dimension of the representation $u_{\alpha}$ as a $k-$vector
space. Furthermore, let us set 
\[
\langle u_{\alpha}\rangle =s^{dim (u_{\alpha})-\epsilon(\alpha)} u_{\alpha}.
\]

For conveniences, we may sometimes simply denote the element $u_{\alpha}$ by 
$\alpha$ for any $\alpha \in \P$, and denote $F_{u_{\alpha}u_{\beta}}^{u_{\lambda}}(rs^{-1})$ by
$g_{\alpha\beta}^{\lambda}$, if no confusion arises. In the rest of this subsection, we will carry out 
all the computations in terms of $\alpha$. It is obvious that the set $\{\langle \alpha \rangle \mid \alpha 
\in \mathcal{P}\}$ is also a $\Q(r,s)-$basis for the algebra $H_{r,s}(A_{\infty})$. Note that we have 
$\langle \alpha_{i} \rangle =\alpha_{i}$ for any given element $\alpha_{i}\in \P$ corresponding to 
the simple root $\alpha_{i}, i\geq 1$. As a result, we can rewrite the multiplication of $H_{r,s}(A_{\infty})$ in 
terms of this new basis as follows
\[
\langle \alpha \rangle \langle \beta
  \rangle=s^{-\epsilon(\alpha)-\epsilon(\beta)-\langle \bf{dim} \alpha,
  \bf{dim} \beta \rangle }\sum _{\lambda \in
  \mathcal{P}}s^{\epsilon(\lambda)}g_{\alpha \beta}^{\lambda}\langle \lambda\rangle
\]
for any $\alpha, \beta \in \mathcal{P}$.
 
In addition, let us denote by 
\[
e(\alpha, \beta)=dim_{k}Hom_{A_{\infty}-mod}(M(\alpha),M(\beta))
\]
and 
\[
 \zeta(\alpha,\beta)= dim_{k}Ext_{A_{\infty}-mod}^{1}(M(\alpha),
M(\beta)).
\]

Recall that Hou and Ye have given an explicit total ordering on the set 
of all isoclasses of finite dimensional indecomposable representations of 
the infinite linear quiver $A_{\infty}$ and used it to construct a PBW base 
for the generic one-parameter Ringel-Hall algebra $H_{q}(A_{\infty})$. 
Following \cite{HY}, we will order all the positive roots as follows:
\[
\a_{11}<\a_{12}<\cdots<\a_{22}<\a_{23}<\cdots.
\]
Obviously, we can see that $Hom(M(\a_{ij}), M(\a_{kl}))\neq 0$ implies 
$\a_{ij}>\a_{kl}$, where $M(\a_{ij}), M(\a_{kl})$ are the indecomposable representations corresponding to the positive roots $\a_{ij}, \a_{kl}$ respectively. For more details about the ordering, we refer the reader to \cite{HY, Rin2}. We should mention that we may write the 
positive roots as $\a_{1}, \a_{2}, \a_{3}, \cdots$ instead.

First of all, we have the following proposition.
\begin{prop}
Let $\alpha_{1},\cdots, \alpha_{t}\in \mathcal{P}$ such that for
$i<j$, we have both $\epsilon(\alpha_{j},\alpha_{i})=0$ and
$\zeta(\alpha_{i},\alpha_{j})=0$. Then 
\[
\langle \bigoplus_{i=1}^{t}\alpha_{i}\rangle =\langle
\alpha_{1}\rangle \cdots \langle \alpha_{t}\rangle.
\] 
\end{prop}
 \qed

\begin{thm}
Let $\alpha, \beta \in \mathcal{P}$ such that $e(\beta, \alpha)=0,
\zeta(\alpha, \beta)=0$. Then we have the following
\[
\langle \beta\rangle \langle \alpha\rangle = r^{\langle
  \alpha,\beta\rangle}s^{-\langle \beta,
  \alpha\rangle }\langle \alpha\rangle \langle \beta\rangle  + \sum _{\gamma \in J(\alpha, \beta)}
  c_{\gamma} \langle \gamma\rangle 
\]
with coefficients $c_{\gamma}$ in $\Z[r^{\pm 1},s^{\pm 1}]$ and $J(\alpha,\beta)$ is
the set of all elements $\lambda \in \mathcal{P}$ which are different from
$\alpha\oplus \beta$ and $g_{\alpha \beta}^{\lambda} \neq 0$.
\end{thm}
 \qed

\begin{prop}
For any given $\alpha \in \mathcal{P}$, we have 
\[
\langle \alpha\rangle =\langle \alpha(\a_{1})\a_{1}\rangle \cdots
\langle \alpha(\a_{m})\a_{m}\rangle.
\]
\end{prop}
\qed

Now let us consider the divided powers of $\langle \a\rangle $ by setting
\[
\langle \a\rangle ^{(t)}=\frac{1}{[t]^{!}_{\epsilon(\a)}}\langle
\a\rangle ^{t}
\]
where $[t]_{\epsilon(\a)}^{!}=\prod_{i=1}^{t}\frac{r^{i\epsilon(\a)}-s^{i\epsilon(\a)}}{r^{\epsilon(\a)}-s^{\epsilon(a)}}$.

Then we have the following lemma.
\begin{lem}
Let $\a$ be a positive root and $t\geq 0$ be an integer. Then we have
the following 
\[
\langle t\a\rangle =\langle \a \rangle ^{(t)}.
\]
\end{lem} 
\qed

For each positive root $\a_{i}$, let us define the following symbol
\[
X_{i}=\langle \a_{i}\rangle.
\]

Then we have the following proposition:
\begin{prop}
Let $\alpha \in \mathcal{P}$ and regard $\alpha$ as a map $\alpha \colon \Phi^{+}\longrightarrow \N_{0}$ with finite support. Let us set $\alpha(i)=\alpha(\a_{i})$, then we have the following 
\[
\langle \alpha\rangle =X_{1}^{(\alpha(1))}\cdots  X_{m}^{(\alpha(m))}=
(\prod_{i=1}^{m}\frac{1}{[\alpha(i)]^{!}_{\epsilon(\a_{i})}}) 
X_{1}^{\alpha(1)}\cdots  X_{m}^{\alpha(m)}.
\]
\end{prop}
\qed

\begin{thm}
The monomials $X_{1}^{\alpha(1)}\cdots X_{m}^{
  \alpha(m)}$ with $\alpha(1), \cdots, \alpha(m)\in \N_{0}$ form a
  $\Q(r,s)-$basis of $H_{r,s}(\Lambda)$; and for $i<j$, we have  
\begin{eqnarray*}
X_{j}X_{i}&= & r^{\langle \dim X_{i}, \dim X_{j}\rangle }s^{-\langle \dim
  X_{j}, \dim X_{i}\rangle } X_{i}X_{j}\\
&& + \sum _{I(i, j)} c(a_{i+1},
  \cdots, a_{j-1}) X_{i+1}^{a_{i+1}}\cdots
  X_{j-1}^{a_{j-1}}
\end{eqnarray*}
with coefficients $c(a_{i+1},\cdots,a_{j-1})$ in $\Q(r,s)$. Here the 
index set $I(i,j)$ is the  set of sequences $(a_{i+1}, \cdots
a_{j-1})$ of natural numbers such that $\sum_{t=i+1}^{j-1}a_{t}{\bf a_{t}=a_{i}+a_{j}}$.
\end{thm}
\qed

Now we define some algebra automorphisms and skew derivations 
on $H_{r,s}(A_{\infty})$. For any $d\in \Z^{\oplus \N}$, we define an algebra 
automorphism $l_{d}$ of $H_{r,s}(A_{\infty})$ as follows
\[
l_{d}(w)=r^{<\dim w, d>}s^{-<d, \dim w>}w
\]
where $w$ is any homogeneous element of
$H_{r,s}(A_{\infty})$. 

Let $H_{j}$ denote the $\Q(r,s)-$subalgebra of $H_{r,s}(A_{\infty})$ generated by the 
generators $X_{1},\cdots, X_{j}$. Thus we have $H_{0}=\Q(r,s)$ and for any $0\leq j\leq m$, 
we have following 
\[
H_{j}=H_{j-1}[X_{j}, l_{j}, \delta_{j}]
\] 
with the automorphism $l_{j}$ and the $l_{j}-$derivation $\delta_{j}$
of $H_{j-1}$. Note that the automorphism $l_{j}$ can be explicitly 
defined as follows
\[
l_{j}(X_{i})=r^{\langle \dim X_{i}, \dim X_{j}\rangle }s^{-\langle \dim X_{j},
  \dim X_{i}\rangle } X_{i}
\]
for $i<j$. And the skew derivation $\delta_{j}$ can be defined as follows:
\[
\delta_{j}(X_{i})=X_{j} X_{i}-l_{j}(X_{i}) X_{j}=\sum _{I(i, j)} c(a_{i+1},
  \cdots, a_{j-1}) X_{i+1}^{a_{i+1}}\cdots
  X_{j_1}^{a_{j-1}}.
\]

It is easy to check that we have the following result.
\begin{prop}
The automorphism $l_{j}$ and the skew derivation $\delta_{j}$ satisfy
the following relation
\[
l_{j}\delta_{j}= r^{\langle \a_{j}, \a_{j}\rangle }s^{-\langle \a_{j},
  \a_{j}\rangle }\delta_{j}l_{j}.
\]
\end{prop}
\qed

\begin{thm}
The two-parameter Ringel-Hall algebra $H_{r,s}(A_{\infty})$ can be
presented as an iterated skew polynomial ring. 
\end{thm}
\qed

\subsection{An algebra isomorphism from $U_{r,s}^{+}(\frak \g)$ 
onto $H_{r,s}(A_{\infty})$}

In this subsection, we are going to establish a graded algebra 
isomorphism from the two-parameter quantized enveloping 
algebra $U_{r,s}^{+}(\mathfrak \g)$ onto the two-parameter 
Ringel-Hall algebra $H_{r,s}(A_{\infty})$. Via this isomorphism, 
all results established in the previous subsection on 
$H_{r,s}(A_{\infty})$ can be transformed to the two-parameter 
quantized enveloping algebra $U_{r,s}^{+}(\mathfrak \g)$. Indeed, 
the isomorphism from $U_{r,s}^{+}(\mathfrak \g)$ onto $H_{r,s}(A_{\infty})$
is the direct limit of the isomorphisms from $U_{r,s}^{+}(\mathfrak sl_{n+1})$ 
onto $H_{r,s}(A_{n})$.

First of all, one can prove the following result, which induces a
homomorphism from $U_{r,s}^{+}(\mathfrak \g)$ into $H_{r,s}(A_{\infty})$.
\begin{lem} Let $\alpha_{i}\in \P$ correspond to the simple 
module $S_{i}$, then we have the following identities in 
$H_{r,s}(\Lambda_{\infty})$.
\begin{eqnarray*}
\alpha_{i}^{2}\alpha_{i+1}^{2}-(r+s)\alpha_{i}\alpha_{i+1}\alpha_{i}+rs\alpha_{i+1}\alpha_{i}^{2}=0,\\
\alpha_{i}\alpha_{i+1}^{2}-(r+s)\alpha_{i+1}\alpha_{i}\alpha_{i+1}+rs\alpha_{i}\alpha_{i+1}^{2}=0,
\end{eqnarray*}
for $i=1, 2, 3,\cdots$.
\end{lem}
{\bf Proof:} Note that we can regard $\alpha_{i}, \alpha_{i+1}$ as elements of the 
two-parameter Ringel-Hall algebra $H_{r,s}(A_{i+1})$, which is a subalgebra of $H_{r,s}(A_{\infty})$. By the result in \cite{T1}, we know that these identities hold in the algebra 
$H_{r,s}(A_{i+1})$. Therefore, we have proved the result as desired.
\qed
 
Now we have the following result which relates Ringel-Hall $H_{r,s}(A_{\infty})$ to the algebra $U_{r,s}^{+}(\mathfrak \g)$.
\begin{thm}
The map 
\[
\eta\colon e_{i} \longrightarrow \alpha_{i}
\]
extends to a $\Q(r,s)-$algebra isomorphism
\[
\eta: U_{r,s}^{+}(\mathfrak \g)\longrightarrow H_{r,s}(A_{\infty}).
\]
\end {thm}
{\bf Proof:} (The proof is essentially borrowed from \cite{Rei} and we include it for completeness. See also \cite{T1}). First of all, note that the quantum Serre relations of $U_{r,s}^{+}(\mathfrak \g)$ are preserved by the map $\eta$. Thus the map $\eta$ does defines an algebra homomorphism from the two-parameter quantized enveloping algebra $U_{r,s}^{+}(\mathfrak \g )$ into the two-parameter twisted Ringel-Hall algebra $H_{r,s}(A_{\infty})$. Now it suffices to show that the map $\eta$ is indeed a bijection.

We first show that the map $\eta$ is surjective by verifying that the algebra $H_{r,s}(A_{\infty})$ is generated by the elements $u_{i}$ which correspond to the irreducible representation $S_{i}$ of the infinite quiver $A_{\infty}$. Let $u_{\alpha}$ be any element in $H_{r,s}(A_{\infty})$, then we can regard $u_{\alpha}$ as an element of a certain subalgebra $H_{r,s}(A_{n})$. Thus we can restrict our proof to the subalgebra $H_{r,s}(A_{n})$. As a result, we have the following:
\[
u_{\alpha}=(\prod_{i=1}^{m}\frac{1}{[\alpha(i)]^{!}_{\epsilon(\a_{i})}}) 
u_{\a_{1}}^{ \alpha(\a_{1})} \cdots u_{\a_{m}}^{ \alpha(\a_{m})}.
\]

Now we need to prove that $u_{\alpha}$ is generated by $u_{i}$ 
for any $\alpha$ corresponding to an indecomposable representations. We prove 
this claim by using induction. Note that $\zeta(\alpha, \alpha)=0$, 
thus we have the following 

\[
u_{\alpha}=u_{1}^{d_{1}}\cdots u_{n}^{d_{n}}-\sum_{\beta\neq \alpha
  \, \dim \beta =\dim \alpha} s^{\langle \beta, \beta\rangle } u_{\beta}. 
\]

However, one sees that the dimension of the module $u_{\beta}$ is 
less than the dimension of the module $u_{\alpha}$. Thus by induction 
on the dimension, we can reduce to the case where $dim
(u_{\alpha})=1$. In this case, the only possibility is that 
$u_{\alpha}=u_{i}$ for some $i$. Thus we have proved the 
statement that every $u_{\alpha}$ is generated by $u_{i}$, 
which further implies that the map $\eta$ is a surjective 
map. We also note that the map $\eta$ is a graded map.

Finally, we show that the map $\eta$ is also injective. Recall that 
$\mathcal{A}=\Q[r,s]_{(r-1,s-1)}$ denote the localization of the polynomial ring 
$\Q[r,s]$ at the maximal ideal $(r-1,s-1)$. Then we know that $\mathcal{A}
=\Q[r,s]_{(r-1,s-1)}$ is a local ring with the residue field $\Q$ and the fractional 
field $\Q(r,s)$. Let $U_{\mathcal{A}}^{+}$ denote the free $\mathcal{A}-$algebra 
generated by the generators $e_{i}$ subject to the quantum Serre relations holding in $U_{r,s}^{+}(\mathfrak \g)$. Also let $U_{\Q}^{+}(\mathfrak \g)$ denote the universal enveloping algebra of the corresponding nilpotent Lie subalgebra $\mathfrak n^{+}$ of $\g$ defined over the base field $\Q$. Then we have the following 
\[
U_{r,s}^{+}(\mathfrak \g)= \Q(r,s)\otimes_{\mathcal{A}} U_{\mathcal{A}}^{+}, 
\quad U_{\Q}^{+}(\mathfrak \g)= \Q \otimes_{\mathcal{A}} U_{\mathcal{A}}^{+}.
\]

For any $\beta \in \Z^{\oplus \N}$, we have the following result via 
Nakayama's Lemma
\begin{eqnarray*}
dim_{\Q} U_{\Q}^{+}(\mathfrak \g)_{\beta} & = & dim_{\Q} (\Q \otimes_{\mathcal{A}}
U_{\mathcal{A}}^{+}))_{\beta}\\
&\geq & dim_{\Q(r, s)}(\Q(r,s)\otimes_{\mathcal{A}} U_{\mathcal{A}}^{+})_{\beta}\\ 
&=& dim_{\Q(r,s)} U_{r,s}^{+}(\mathfrak \g)_{\beta}\\
&\geq & dim_{\Q(r,s)}H_{r,s}(A_{\infty})_{\beta}.
\end{eqnarray*}

Note that we also have the following result: 
\[
dim_{\Q} U_{\Q}^{+}(\mathfrak \g)_{\beta}=dim_{\Q(r,s)}H_{r,s}(A_{\infty})_{\beta}.
\]

Thus we have proved that the map $\eta$ is injective. Therefore, 
the map $\eta$ is an algebra isomorphism from $U_{r,s}^{+}(\mathfrak \g)$ 
onto $H_{r,s}(A_{\infty})$ as desired.
\qed
 
\qed
 
Based on the previous theorem, the following corollary is in order.
\begin{cor}
The algebra $U_{r,s}^{+}(\mathfrak \g)$ has a $\Q(r,s)-$basis 
parameterized by the isomorphism classes of all finite dimensional 
representations of the infinite quiver $A_{\infty}$. In particular, we have 
\[
U_{r,s}^{+}(\mathfrak sl_{\infty})=\lim_{n\mapsto \infty} U_{r,s}^{+}(\mathfrak sl_{n+1}).
\]
\end{cor}
\qed

\section{The extended two-parameter Ringel-Hall algebras 
$\overline{H_{r,s}(A_{\infty})}$}

For the purpose of realizing the Borel subalgebra $U_{r,s}^{\geq
  0}(\mathfrak \g)$ of the two-parameter quantum group
$U_{r,s}(\mathfrak \g)$, we define the extended Ringel-Hall algebra 
$\overline{H_{r,s}(A_{\infty})}$ by adding the torus part. 
In particular, we show that there is a Hopf algebra structure on 
this extended two-parameter Ringel-Hall algebra
  $\overline{H_{r,s}(A_{\infty})}$; as a result we prove that 
$U_{r,s}^{\geq 0}(\mathfrak \g)$ is isomorphic to the extended 
two-parameter Ringel-Hall algebra $\overline{H_{r,s}(A_{\infty})}$ 
as a Hopf algebra. Similarly, we can use an extended two-parameter 
Ringel-Hall algebra to realize the Borel subalgebra 
$U^{\leq 0}(\mathfrak sl_{\infty})$. Therefore, we will obtain a PBW-basis of 
two-parameter quantum group $U_{r,s}(\mathfrak \g)$.

\subsection{Extended Ringel-Hall algebras $\overline{H_{r,s}(A_{\infty})}$}
Let us define $\overline{H_{r,s}(A_{\infty})}$ to be a free $\Q(r,s)-$module 
with the following basis
\[
\{k_{\alpha}u_{\lambda}\mid \alpha \in \Z[I], \, \lambda \in \P \}.
\]

Moreover one will define an algebra structure on the module 
$\overline{H_{r,s}(A_{\infty})}$ as follows.
\begin{eqnarray*}
u_{\alpha}u_{\beta}=\sum_{\lambda \in \P}
s^{-\langle \alpha,\beta \rangle }F^{u_{\lambda}}_{u_{\alpha}, u_{\beta}}(rs^{-1})u_{\lambda},\quad \text{for
  any}\, \alpha, \beta \in \P,\\
k_{\alpha}u_{\beta}=r^{\langle \beta,
  \alpha\rangle }s^{-\langle \alpha,\beta\rangle }u_{\beta}k_{\alpha} \quad \text{for any
}\quad \alpha \in \Z[I],\beta \in \P,\\
k_{\alpha}k_{\beta}=k_{\beta}k_{\alpha} \quad \text{for any}
  \quad \alpha, \beta \in \Z[I].
\end{eqnarray*}

Indeed, we have the following
\begin{prop} 
For any elements $x, y, z \in \Z[I]$ and $\alpha , \beta, \gamma 
\in \mathcal{P}$, we have the following
\[
[(k_{x}u_{\alpha})(k_{y}u_{\beta})](k_{x}u_{\alpha})=(k_{x}u_{\alpha})[(k_{y}u_{\beta})(k_{z}u_{\gamma})].
\]
In particular, with the above defined multiplication,
$\overline{H_{r,s}(\Lambda_{\infty})}$ is an associative $\Q(r,s)-$algebra.
\end{prop}
{\bf Proof:} Once we choose $x,y,z,$ and $\alpha,\beta,\gamma$, we can
restrict to the subgroup $\overline{H_{r,s}(A_{m})}$ of
$\overline{H_{r,s}(A_{\infty})}$ for some $m$. Since $\overline{H_{r,s}(A_{m})}$ is 
an associative algebra with the restricted multiplication, thus we have proved all the statements.
\qed

Furthermore, we have the following result.
\begin{thm}
The map $\eta$ extends to a $\Q(r,s)-$algebra isomorphism from 
$U_{r,s}^{\geq 0}(\mathfrak \g)$ onto $\overline{H_{r,s}(\Lambda_{\infty})}$ 
via the map $\eta(w_{i})=k_{i}$ and $\eta(e_{i})=u_{\alpha_{i}}$.
\end{thm}
{\bf Proof:} The proof is straightforward.
\qed

As a result, we have the following description about a basis for 
the algebra $U_{r,s}^{+}(\frak \g)$
\begin{cor}
The set ${\bf B}^{+}=\{w_{\alpha}\eta^{-1}(u_{\lambda})\mid \alpha \in \Z[\N], \lambda \in
\P \}
$ is a $\Q(r,s)-$basis of $U_{r,s}^{\geq 0}(\frak \g)$.
\end{cor}
\qed

\subsection{A Hopf algebra structure on $\overline{H_{r,s}(A_{\infty})}$}
Now we are going to introduce a Hopf algebra structure on the extended
two-parameter Ringel-Hall algebra $\overline{H_{r,s}(A_{\infty})}$. In 
particular, we have the following result.
\begin{thm}
The algebra $\overline{H_{r,s}(A_{\infty})}$ is a Hopf algebra with the 
Hopf algebra structure defined as follows.
\begin{enumerate}
\item Multiplication:
\begin{eqnarray*}
u_{\alpha}u_{\beta}=\sum_{\lambda \in \P}
s^{-\langle \alpha,\beta \rangle }g^{\lambda}_{\alpha \beta}u_{\lambda}\quad \text{for
  any}\quad \alpha, \beta \in \B,\\
k_{\alpha}u_{\beta}=r^{\langle \beta,
  \alpha \rangle }s^{-\langle \alpha,\beta\rangle }u_{\beta}k_{\alpha} \quad \text{for any
}\quad \alpha \in \Z[I],\beta \in \P,\\
k_{\alpha}k_{\beta}=k_{\beta}k_{\alpha} \quad \text{for any}
  \quad \alpha, \beta \in \Z[I].
\end{eqnarray*}

\item  Comultiplication:
\begin{eqnarray*}
\Delta(u_{\lambda})=\sum_{\alpha,\beta \in \P}r^{\langle \alpha,
    \beta\rangle }
  \frac{a_{\alpha}a_{\beta}}{a_{\lambda}}g^{\lambda}_{\alpha
\beta}u_{\alpha}k_{\beta}\otimes u_{\beta}\quad \text{for any } \quad \lambda \in \P,\\
\Delta(k_{\alpha})=k_{\alpha}\otimes k_{\alpha} \quad \text{for
    any} \quad \alpha \in \Z[I].
\end{eqnarray*}
\item Counit:
\begin{eqnarray*}
\epsilon(u_{\lambda})=0 \quad \text{for all} \quad \lambda \neq 0 \quad
\text{and} \quad \epsilon(k_{\alpha})=1 \quad \text{for any} \quad \alpha \in \P.
\end{eqnarray*}
\item Antipode:

\begin{eqnarray*}
\sigma(u_{\lambda})&=&\delta_{\lambda,0}+\sum_{m\geq 1}(-1)^{m}
\times\sum_{\pi \in \P, \lambda_{1},\lambda_{2},\cdots,\lambda_{m}\in \P_{1}}
(rs^{-1})^{\sum_{i<j}\langle \lambda_{i},\lambda_{j}\rangle}\\
& & \frac{a_{\lambda_{1}}\cdots
    a_{\lambda_{m}}}{a_{\lambda}} g^{\lambda}_{\lambda_{1} \cdots
  \lambda_{m}} g^{\pi}_{\lambda_{1}\cdots\lambda
  _{m}}k_{-\lambda}u_{\pi}
\end{eqnarray*}
for any  element $\lambda \in \P$ and 

\[\sigma(k_{\alpha})=k_{-\alpha} \quad \text{for any} \quad \alpha \in
\Z[I].\\
\]
\end{enumerate}
In particular, we have the following
\[
\overline{H_{r,s}(A_{\infty})}=\lim_{n\mapsto \infty} \overline{H_{r,s}(A_{n})}
\]
as a direct limit of Hopf subalgebras.
\end{thm}
\qed

The proof of the above theorem consists of a couple of lemmas which can be proved as the 
finite dimensional case. And we refer the reader to \cite{T1, Xiao} for more details. Namely, we have the following lemmas.
\begin{lem}
The comultiplication $\Delta$ is an algebra endomorphism of 
$\overline{H_{r,s}(A_{\infty})}$.
\end{lem}
\qed
 
\begin{lem}
For any $\lambda \in \mathcal{P}$, we have the following 
\[
\mu (\sigma \otimes 1) \Delta(u_{\lambda}) = \delta_{\lambda 0}
\]
and 
\[
\mu (1\otimes \sigma)\Delta(u_{\lambda})=\delta_{\lambda 0}.
\] 
\end{lem}
\qed
 
\subsection{A Hopf algebra isomorphism from $U_{r,s}^{\geq 0}(\mathfrak
  \g)$ onto $\overline{H_{r,s}(A_{\infty})}$}
In this subsection, we will prove that the Borel subalgebras 
$U_{r,s}^{\geq 0}(\mathfrak \g)$ and $U_{r,s}^{\leq 0}(\mathfrak \g)$ of the 
two-parameter quantum group $U_{r,s}(\mathfrak \g)$ can be realized as the 
extended two-parameter Ringel-Hall algebra $\overline{H_{r,s}(A_{\infty})}$ 
and $\overline{H_{s^{-1},r^{-1}}(A_{\infty})}$ as Hopf algebras. As a result, 
we shall derive a PBW-basis for the two-parameter quantum group  $U_{r,s}(\mathfrak \g)$.

\begin{thm}
We have that 
\[
U_{r,s}^{\geq 0}(\mathfrak \g)\cong \overline{H_{r,s}(A_{\infty})}
\]
and
\[
U_{r,s}^{\leq 0}(\mathfrak \g)\cong
\overline{H_{s^{-1},r^{-1}}(A_{\infty})}
\]
as Hopf algebras.
\end{thm}
\qed

Let ${\bf B}^{-}$ denote the $\Q(r,s)-$basis constructed for the algebra 
$U_{r,s}^{\leq 0}(\mathfrak \g)$ via the Ringel-Hall algebra $\overline{H_{s^{-1},r^{-1}}(A_{\infty})}$, then we have the following:
\begin{cor}
The set ${\bf B}^{+}\times {\bf B}^{-}$ is a $\Q(r,s)-$basis for the
two-parameter quantum groups $U_{r,s}(\mathfrak \g)$.
\end{cor}
\qed

Furthermore, we have the following result, which provides a 
bridge from the finite dimensional case to the infinite case. 
\begin{thm}
The two-parameter quantum group $U_{r,s}(\mathfrak \g)$ is the direct 
limit of the direct system $\{U_{r,s}(\mathfrak sl_{n+1})\mid n\in
\N\}$ of the Hopf subalgebras $U_{r,s}(\mathfrak sl_{n+1})$ of $U_{r,s}(\mathfrak \g)$ . That is
\[
U_{r,s}(\mathfrak \g)=\lim_{n\mapsto \infty}U_{r,s}(\mathfrak sl_{n+1}).
\]
In particular, we have
\begin{eqnarray*}
U^{\pm 1}_{r,s}(\mathfrak \g)=\lim_{n\mapsto \infty}U^{\pm 1}_{r,s}(\mathfrak
sl_{n+1}),\\
U^{0}_{r,s}(\mathfrak \g)=\lim_{n\mapsto \infty}U^{0}_{r,s}(\mathfrak
sl_{n+1}),\\
U^{\geq 0}_{r,s}(\mathfrak \g)=\lim_{n\mapsto \infty}U^{\geq 0}_{r,s}(\mathfrak sl_{n+1}),\\
U^{\leq 0}_{r,s}(\mathfrak \g)=\lim_{n\mapsto \infty }U^{\leq
  0}_{r,s}(\mathfrak sl_{n+1}).
\end{eqnarray*}
\end{thm}
{\bf Proof:} It is obvious that $U_{r,s}(\mathfrak sl_{n+1})$ are
Hopf subalgebras of $U_{r,s}(\mathfrak sl_{\infty})$ and $U_{r,s}(\mathfrak
sl_{m+1})$ for $m\geq n$. In addition, any element of
$U_{r,s}(\mathfrak sl_{\infty})$ is an element of a certain 
$U_{r,s}(\mathfrak sl_{n+1})$. Thus we are done with the proof.
\qed

\section{Monomial bases and bar-invariant bases of $U_{r,s}^{+}(\mathfrak sl_{\infty})$} In this section, we study various bases of $U_{r,s}^{+}(\mathfrak sl_{\infty})$  via the theory of generic extensions. Note that the construction of monomial bases using generic extension theory for the Ringel-Hall algebras of type $A, D, E$ has been done in \cite{DD1}. The idea of the construction 
is to use the monoidal structure on the set $\mathcal{M}$ of isoclasses of finite 
dimensional representations of the corresponding quiver $\mathcal{Q}$ and the 
Bruhat-Chevalley type partial ordering in $\mathcal{M}$. Note that the arguments used in \cite{DD1} can be completely transformed to the case of $\mathfrak sl_{\infty}$. Therefore, 
we will state most of the results for monomial bases without much detail. For the details, we 
refer the reader to \cite{DD1, Rei}. 

For the reader's convenience, we will recall the necessary details about the 
the generic extensions from \cite{DD1, Rei}. Note that there exists a bijective 
correspondence between the set of positive roots $\Phi^{+}$ of the root system $\Phi$ 
associated to $\mathfrak sl_{\infty}$ and the set of isoclasses of finite dimensional 
indecomposable representations of $A_{\infty}$. For any $\beta \in \Phi^{+}$, 
we will denote by $M(\beta)=M_{k}(\beta)$ the corresponding indecomposable 
representation of $A_{\infty}$. By the Krull--Remak-Schmidt theorem, we shall 
have the following
\[
M(\lambda)=M_{k}(\lambda)\colon = \bigoplus _{\beta \in \Phi^{+}}\lambda(\beta)M_{k}(\beta)
\]
for some function $\lambda \colon \Phi^{+}\longrightarrow
\mathbb{N}_{0}$ with a finite support. Therefore, the isoclasses of finite dimensional representations of $A_{\infty}$ are indexed by the following set 
\[
\Lambda =\{\lambda \colon \Phi^{+}\longrightarrow \mathbb{N}\,\text{with a finite support}\}\cong \mathbb{N}_{0}^{\oplus \Phi^{+}}.
\]
From now on, we will use the set $\Lambda$ to index the objects of the category $A_{\infty}-$mod. 

Next, we are going to recall some information about generic extensions of representations of Dynkin quivers. We should mention that all the arguments used in the finite dimensional cases 
of type $A, D, E$ can be transformed to the $\mathfrak sl_{\infty}$. We refer the interested reader to the references \cite{DD1, Rei} for details.

Let us fix $k$ to algebraically closed. Let us denote by $\mathcal{Q}=(\mathcal{Q}_{0}, \mathcal{Q}_{1})$ the quiver $A_{\infty}$. Fix a ${\bf d}=(d_{i})_{i}\in \N_{0}^{\oplus \Phi^{+}}$ and we may choose $n$ large enough so that ${\bf d}$ can be regarded as an 
element in $\N_{0}^{n}$. For any given ${\bf d}$, we can define an affine space as follows
\[
R({\bf d})=R(\mathcal{Q},{\bf d})\colon=\prod_{\alpha \in
    \mathcal{Q}_{1}} Hom_{k}(k^{d_{t(\alpha)}},
    k^{d_{h_{\alpha}}})\cong\prod_{\alpha \in
    \mathcal{Q}_{1}} {k^{d_{t_{\alpha}}\times d_{h_{\alpha}}}}.
\] 
Thus, a point $x=(x_{\alpha})_{\alpha}$ of $R({\bf d})$ determines a
finite dimensional representation $V(x)$ of $\mathcal{Q}=A_{\infty}$. 
The algebraic group $GL({\bf d})=\prod_{i=1}^{n}GL_{d_{i}}(k)$ acts on the space $R({\bf d})$ by the conjugation
\[
(g_{i})_{i}\dots(x_{\alpha})_{\alpha}=(g_{h(\alpha)})x_{\alpha}g_{t(\alpha)}^{-1})_{\alpha}.
\]
and the $GL({\bf d})-$orbits $\mathcal{O}_{x}$ in $R({\bf d})$ correspond
bijectively to the isoclasses $[V(x)]$ of finite dimensional representations of
$\mathcal{Q}$ with the dimension vector ${\bf d}$. The stabilizer $GL({\bf
d})_{x}=\{g\in GL({\bf d})\mid gx=x\}$ of $x$ is the group of
automorphisms of $M\colon =V(x)$ which is zariski-open in
$End_{A_{n}-mod}(M)$ and has a dimension equal to the
$dim_{A_{n}-mod}(M)$. It follows that the orbit
$\mathcal{O}_{M}\colon = \mathcal{O}_{x}$ of $M$ has a dimension 
\[
dim\mathcal{O}_{M}=dimGL({\bf d})-dimEnd_{A_{n}-mod}(M).
\]

Now we have the following result, whose proof is the same as the one in \cite{Rei}.
\begin{lem}
For $x\in R({\bf d_{1}})$ and any $y\in R({\bf d_{2}})$, let
$\mathcal{E}(\mathcal{O}_{x},\mathcal{O}_{y})$ be the set of all $z\in
R({\bf d})$ where ${\bf d=d_{1}+d_{2}}$ such that $V(z)$ is an
extension of some $M\in \mathcal{O}_{x}$ by some $N\in
\mathcal{O}_{y}$. Then $\mathcal{E}(\mathcal{O}_{x},\mathcal{O}_{y})$
is irreducible. 
\end{lem}
\qed

Given any two finite-dimensional representations of $M, N$ of the infinite linear 
quiver $A_{\infty}$, let us consider the extensions 
\[
0\longrightarrow N\longrightarrow L\longrightarrow M\longrightarrow 0
\]
of $M$ by $N$. By the lemma, there is a unique (up to isomorphism)
such extension $G$ with $dim\mathcal{O}_{G}$ being maximal. We call
$G$ the generic extension of $M$ by $N$, and denoted by $M\ast N$.
For any two representations $M,N$, we say $M$ degenerates to $N$, or
that $N$ is a degeneration of $M$, and write  $[N]\leq [M]$ (or simply
$N\leq M$) if $\mathcal{O}_{N}\subseteq \overline{\mathcal{O}_{M}}$ 
which is the closure of $\mathcal{O}_{M}$. Note that $N< M$ if and
only if $\mathcal{O}_{N}\subset \overline{\mathcal{O}_{M}}\backslash
\mathcal{O}_{M}$.

Similar to the result in \cite{DD1, Rei}, one knows that the relation $\leq $ is 
independent of the base field $k$ and it provides a partial order on the set 
$\Lambda$ via setting $\lambda \leq \mu $ if and only if $M_{k}(\lambda)\leq M_{k}(\mu)$ for 
any given algebraically closed field $k$.

Using the same arguments as in \cite{DD1, Rei},  we shall have the following result.
\begin{thm} (1). If $0\longrightarrow N\longrightarrow E\longrightarrow M\longrightarrow 0$ is  exact and non-split, then $M\oplus N<E$.\\
(2). Let $M, N, X$ be finite dimensional representations of the quiver $A_{\infty}$. Then  $X\leq M\ast N$ if and only if there exit $M^{\prime}\leq M, N^{\prime}\leq N$ such that $X$ is an extension of $M^{\prime}$ by $N^{\prime}$. In particular, we have $M^{\prime}\leq M, N^{\prime}\leq N \Longrightarrow M^{\prime}\ast N^{\prime} \leq M\ast N$.\\
(3). Let $\mathcal{M}$ be the set of isoclasses of finite dimensional representations of $A_{\infty}$ and define a multiplication $\ast$ on $\mathcal{M}$ by $[M]\ast[N]=[M\ast N]$ for any $[M], [N]\in \mathcal{M}$. Then $\mathcal{M}$ is  a monoid with identity $1=[0]$ and the multiplication $\ast$ preserves the induced partial ordering on $\mathcal{M}$. \\
(4). $\mathcal{M}$ is generated by irreducible representations $[S_{i}], i\in I$ subject to the following relations
\begin{enumerate}
\item $[E_{i}]\ast[E_{j}]=[E_{j}][E_{i}]$ if $i, j$ are not connected by an arrow,
\item $[E_{i}]\ast[E_{j}]\ast[E_{i}]=[E_{i}]\ast[E_{i}]\ast[E_{j}]$ and $[E_{j}]\ast[E_{i}]\ast[E_{j}]=[E_{i}]\ast[E_{j}]\ast[E_{j}]$ if there exists an arrow from $i$ to $j$.
\end{enumerate}
\end{thm}
\qed

Let us denote by $\Omega$ the set of all words formed by letters in $I$. It is easy to see that 
for any given word $w=w_{1}\cdots w_{m}\in \Omega$, we can set the following finite dimensional representations of $A_{\infty}$
\[
M(w)=S_{w_{1}}\ast S_{w_{2}}\ast \cdots \ast S_{w_{m}}.
\]

Note that there is a unique $M(\mathfrak{p}(w)) \in A_{\infty}-mod$ such that $M(w)\cong M(\mathfrak{p}(w))$, which enables us 
to define a function as follows
\[
\mathfrak{p}\colon \Omega \longrightarrow A_{\infty}-mod, w\mapsto M(\mathfrak{p}(w)).
\]

Furthermore, we shall have the following  result on this function.
\begin{thm}
The map $\mathfrak{p}$ induces a surjection 
\[
\mathfrak{p}\colon \Omega \longrightarrow A_{\infty}-mod.
\] 
\end{thm}
{\bf Proof:} Once again, we can restrict the function to a certain subcategory $A_{m}-$mod, where the property holds.
\qed

Therefore, $\mathfrak{p}$ induces a partition of the set $\Omega=\cup_{\lambda \in \Lambda} \Omega_{\lambda}$ with $\Omega_{\lambda}=\mathfrak{p}^{-}(\lambda)$. We will call each $\Omega_{\lambda}$ a fiber of the map $\mathfrak p$. 

Now we are going to recall some information on the partial ordering $\leq$. Let 
$w=i_{1}\cdots i_{m}$ be a word in $\Omega$. Then $w$ can be uniquely expressed in the tight form $w=j_{1}^{e_{1}}\cdots
j_{t}^{e_{t}}$ where $e_{r}\geq 1, 1\leq r\leq t$, and $j_{r}\neq
j_{r+1}$ for $1\leq r\leq t-1$. A filtration
\[
0=M_{t}\subset M_{t-1}\subset \cdots M_{1}\subset M_{0}=M
\]
of a module $M$ is called a reduced filtration of type $w$ if
$M_{r-1}/M_{r}\cong e_{r}S_{r}$, for all $1\leq r\leq t$. Note that
any reduced filtration of $M$ of type $w$ can be refined to a
composition of $M$ of type $w$. Conversely, given a composition series
of $M$, there is a unique reduced filtration of $M$. Let us denote by
$\varphi_{w}^{\lambda}(x)$ the Hall polynomial
$\varphi_{\mu_{1}\cdots\mu_{t}}^{\lambda}(x)$ where
$M(\mu_{r})=e_{r}S_{r}$. Let us denote by $\gamma_{w}^{\lambda}(q_{k})$ 
the number of the reduced filtrations of $M_{k}(\lambda)$ over the base 
field $k$ when $k$ is a finite field. A word $w$ is called distinguished if $\gamma_{w}^{\mathfrak{p}(w)}=1$. Note that $w$ is distinguished if 
and only if, for some algebraically closed field $k$, $M_{k}(\mathfrak{p}(w))$ 
has a unique reduced filtration of type $w$. Similar to \cite{DD1}, we have the 
following results.
\begin{lem}
(See also {\bf Lemma 4.1} in \cite{DD1}) Let $\omega \in \Omega$ and $\mu \geq \lambda $ in $\Lambda$. Then 
$\varphi_{\omega}^{\mu}\neq 0$ implies $\varphi_{\omega}^{\lambda}\neq 0$.
\end{lem}
\qed

\begin{thm}
(See also {\bf Theorem 4.2} in \cite{DD1}) Let $\lambda, \mu \in \Lambda $. Then $\lambda \leq \mu $ if and only
if there exists a word $\omega \in \mathfrak{p}^{-1}(\mu)$ such that 
$\varphi^{\lambda}_{\mu}\neq 0$.
\end{thm}
\qed

\begin{lem}
(See also {\bf Lemma 5.2} in \cite{DD1}) Every fiber of $\mathfrak{p}$ contains a distinguished word.
\end{lem}
\qed

Let us define $[[e_{a}]]^{!}=[[1]]\cdots[[e_{a}]]$ with $[[m]]=\frac{1-(rs^{-1})^{m}}{1-rs^{-1}}$. Then we shall have the following result.
\begin{lem}
(See also {\bf Lemma 6.1} in \cite{DD1}) Let $w\in \Omega$ be a word with the tight form $j_{1}^{e_{1}}\cdots
j_{t}^{e_{t}}$. Then, for each $\lambda \in \Lambda$, 
\[
\varphi_{w}^{\lambda}(rs^{-1})=\gamma_{w}^{\lambda}(rs^{-1})\prod_{r=1}^{t}[[e_{r}]]^{!}.
\]
In particular,
$\varphi_{w}^{\mathfrak{p}(w)}=\prod_{r=1}^{t}[[e_{r}]]^{!}$ if $w$ is distinguished.
\end{lem}
\qed

For any given word $w=i_{1}\cdots c_{m} \in \Omega$, we can associate a monomial
\[
u_{w}=u_{i_{1}}\cdots u_{i_{m}}\in H_{r,s}(A_{\infty}).
\] 

\begin{prop}
For any $w \in \Omega$ with the tight form $j_{1}^{r_{1}}\cdots
j_{t}^{e_{t}}$, we have 
\[
u_{w}=\sum _{\lambda \leq \mathfrak{p}(w)}
\varphi_{w}^{\lambda}(rs^{-1})u_{\lambda}=\prod_{r=1}^{t}[[e_{r}]]^{!}\sum_{\lambda
\leq \mathfrak{p}(w)}\gamma_{w}^{\lambda}(rs^{-1})u_{\lambda}.
\]
Moreover, the coefficients appearing in the sum are all nonzero.
\end{prop}
\qed

As a result, we shall have the following theorem.
\begin{thm}
For each given $\lambda \in \Lambda$, let us choose an arbitrary word
$w_{\lambda}\in \mathfrak{p}^{-1}(w)$. Then the set $\{ u_{w_{\lambda}}\mid
\lambda \in \Lambda\}$ is a $\Q(r,s)-$basis of
$H_{r,s}(A_{\infty})$. Moreover, if all the words are chosen to be distinguished, then this set is a $\Z[r, s]_{(r-1, s-1)}-$basis of $H_{r,s}(A_{\infty})_{\Z[r, s]_{(r-1,s-1)}}$. 
\end{thm}
\qed
 
\subsection{A bar-invariant basis of $U^{+}_{r,s}(\mathfrak sl_{\infty})$} It is easy to see that the algebra $U_{r,s}^{+}(\mathfrak{sl_{\infty}})$ admits a $\Q-$linear involution defined as follows
\[
\overline{r}=s, \overline{s}=r, \overline{e_{i}}=e_{i}\, \text{for all}\, i \in I.
\]
And we will refer this involution as the bar-involution. In this subsection, we will construct a bar-invariant basis for $U_{r,s}^{+}(\mathfrak{sl_{\infty}})$. 

Denote by $[M,N]=dim_{k}Hom(M, N)$ and $[M, N]^{1}=dim_{k}Ext^{1}(M, N)$. Let us set $c_{M, N}^{X}=s^{[X, X]-[M,N]+[M,N]^{1}-[M,M]-[N,N]}F_{M,N}^{X}(rs^{-1})$. It is obvious that the same proof in \cite{Rei} shall yield the following result.
\begin{prop}
Let us write $\overline{u_{\alpha}}=\sum_{\beta}\omega_{\beta}^{\alpha}u_{\beta}$, then we have
\begin{enumerate}
\item $\omega_{\beta}^{\alpha}=0$ unless $\beta \leq \alpha$, and $\omega_{\alpha}^{\alpha}=1$,
\item if $u_{\alpha}=M\oplus N$ for finite dimensional representations $M, N$ with $[N, M]=0=[M, N]^{1}$, then 
\[
\omega_{\beta}^{\alpha}=\sum_{M^{\prime}\leq M, N^{\prime}\leq N} \omega_{M^{\prime}}^{M}\omega_{N^{\prime}}^{N}c_{M^{\prime}N^{\prime}}^{\alpha},
\]
\item if $u_{\alpha}$ is an exponent of a finite dimensional indecomposable representation, then 
\[
\omega_{\beta}^{\alpha}=s^{[u_{\beta},u_{\beta}]^{1}}-\sum_{\beta \leq \gamma < \alpha}r^{[u_{\gamma},u_{\gamma}]^{1}}\omega_{\beta}^{\gamma},
\]
\item $\omega_{\beta}^{\alpha} \in s^{[u_{\beta},u_{\beta}]-[u_{\alpha},u_{\alpha}]}\Z[rs^{-1}]$.
\end{enumerate}
\end{prop}
\qed

Furthermore, using the arguments in \cite{Lus2, Rei}, we shall have the following result on a bar-invariant basis of the algebra $U^{+}_{r,s}(\mathfrak sl_{\infty})$. 
\begin{thm}
For each isoclass $\alpha$, there exists a unique element
\[
\mathcal{C}_{\alpha} \in u_{\alpha}+s^{-1}\Z[rs, r^{-1}s^{-1}, s][{\bf B}\backslash \{u_{\alpha}\})
\]
such that $\overline{\mathcal{C}_{\alpha}}=\mathcal{C}_{\alpha}$. Write $\mathcal{C}_{\alpha}=\sum_{\beta} \zeta^{\alpha}_{\beta}u_{\beta}$, we have 
\begin{enumerate}
\item $\zeta^{\alpha}_{\beta}=0$ unless  $\beta \leq \alpha$, and $\zeta_{\alpha}^{\alpha}=1$,
\item $\zeta_{\beta}^{\alpha} \in s^{[u_{\beta},u_{\beta}]-[u_{\alpha},u_{\alpha}]}\Z[rs^{-1}]$,
\item Denote by $\hat{\zeta}_{\beta}^{\alpha}(v) \in \Z[v, v^{-1}]$ the specialization of $\zeta^{\alpha}_{\beta}$ to $\alpha=v=s^{-1}$, we have 
\[
\zeta_{\beta}^{\alpha}=(\sqrt{rs})^{[u_{\beta},u_{\beta}]-[u_{\alpha},u_{\alpha}]}\hat{\zeta}_{Y}^{X}(\sqrt{rs^{-1}}).
\]
\end{enumerate}
\end{thm}
\qed

{\bf Acknowledgement:} The author would like to thank Sarah Witherspoon for some helpful discussions about this work.

\end{document}